\documentclass[a4paper,11pt]{amsart}

\usepackage[cp1250]{inputenc}
\usepackage{amsmath}                            
\usepackage{amssymb}
\usepackage{amsfonts}
\usepackage{mathrsfs}
\usepackage{cite}
\usepackage{textcomp}
\usepackage{array}
\usepackage{multirow}
\usepackage{float}
\usepackage{graphicx}
\usepackage{setspace}
\usepackage{caption} 
\usepackage{graphicx}
\usepackage{mathptmx}
\usepackage[all]{xy}
\usepackage{fancybox}
\usepackage{footmisc}
\usepackage[most]{tcolorbox}

\newtcolorbox{mybox}[2][]
{colbacktitle=white, coltitle=black, boxed title style={boxrule=0.2mm},
attach boxed title to top center={yshift=-8pt}, 
title = #2,#1,enhanced, boxrule=0.2mm, colframe=black}

\DeclareMathAlphabet{\mathpzc}{OT1}{pzc}{m}{it}

\newtheorem{theorem}{Theorem}[section]
\newtheorem{lemma}[theorem]{Lemma}

\newtheorem{property}[theorem]{}
\theoremstyle{definition}

\theoremstyle{remark}
\newtheorem{remark}[theorem]{Remark}
\newtheorem*{rem}{Remark}

\numberwithin{equation}{section}

\allowdisplaybreaks[4]

\newcolumntype{M}[1]{>{\raggedright}m{#1}}

\begin{document}

\title{On the Zariski multiplicity conjecture for\\ weighted homogeneous and Newton non-degenerate\\ line singularities} 

\author{Christophe Eyral and Maria Aparecida Soares Ruas}

\subjclass[2010]{14B05, 14B07, 14J70, 14J17, 32S15, 32S25, 32S05.}

\keywords{Line singularity; Weighted homogeneity; Newton non-degeneracy; topological equisingularity; equimultiplicity; Iomdine-L\^e-Massey formula.}

\begin{abstract}
We present new families of weighted homogeneous and Newton non-degenerate line singularities that satisfy the Zariski multiplicity conjecture.
\end{abstract}

\maketitle

\markboth{C. Eyral and M. Ruas}{On the Zariski multiplicity conjecture} 

\section{Introduction}\label{intro}

Let $\mathbf{z}:=(z_1,\ldots,z_n)$ be linear coordinates for $\mathbb{C}^n$, and let $f\colon (\mathbb{C} \times \mathbb{C}^n,\mathbb{C} \times \{\mathbf{0}\}) \rightarrow (\mathbb{C},0)$, $(t,\mathbf{z})\mapsto f(t,\mathbf{z})$,
be a polynomial function such that for each $t$ the function  $\mathbf{z}\mapsto f_t(\mathbf{z}):= f(t,\mathbf{z})$ is reduced at~$\mathbf{0}\in\mathbb{C}^n$. 
As usual, we denote by $V(f_t)$ the corresponding hypersurface $f_t^{-1}(0)$ in $\mathbb{C}^n$.
We say that the family $\{f_t\}$ is topologically equisingular (at $\mathbf{0}$)\footnote{In this paper, we are only interested in the behaviour of functions (or hypersurfaces) near the origin $\mathbf{0}\in\mathbb{C}^n$, unless otherwise stated. Hereafter, we shall omit the words ``at $\mathbf{0}$.''} if the ambient topological type of the hypersurface-germ $(V(f_t),\mathbf{0})$ is independent of $t$ for all small $t$, that is, if there exists a family $\{\varphi_t\}$ of germs of self-homeomorphisms of $(\mathbb{C}^n,\mathbf{0})$  depending continuously on $t$ and such that for any $t$ sufficiently small $\varphi_t(V(f_0),\mathbf{0})=(V(f_t),\mathbf{0})$.
The \emph{Zariski multiplicity conjecture} \cite{Z} says that if the family $\{f_t\}$ is topologically equisingular, then it is equimultiple (i.e., the multiplicity (equivalently, the order) of $f_t$ at $\mathbf{0}$ is independent of $t$ for all small $t$). This conjecture, posed almost fifty years ago, has been solved only in a few special cases. For a survey on that question and related topics, we refer the reader to \cite{E2,EyBook}.

For example, among known results, there is the following theorem due to Greuel \cite{G} and O'Shea \cite{O'Sh}.

\begin{theorem}[Greuel and O'Shea]\label{thm-GOSh}
Suppose that $\{f_t\}$ is a family of isolated singularities (i.e., for all small $t$, $f_t$ has an isolated singularity at $\mathbf{0}$) such that the polynomial function $f_0$ is weighted homogeneous with respect to a system of positive integer weights $(w_1,\ldots,w_n)$. If furthermore the Milnor number $\mu(f_t)$ of $f_t$ at~$\mathbf{0}$ is independent of $t$ for all small $t$,\footnote{It is well known that if the family $\{f_t\}$ is topologically equisingular, then the Milnor number $\mu(f_t)$ is independent of $t$ for all small $t$. Both conditions are equivalent if $n\not=3$ (cf.~\cite{LR}).} then $\{f_t\}$ is equimultiple.
\end{theorem}

The proof of Theorem \ref{thm-GOSh} uses a crucial result of Varchenko~\cite{V}.
The case where $f_0$ is homogeneous was first proved by Gabri\`elov and Ku\v{s}nirenko \cite{GK}.

Another well known class of isolated hypersurface singularities which satisfies the Zariski multiplicity conjecture is described by the following theorem due to Abderrahmane \cite{A} and Saia and Tomazella \cite{ST}.

\begin{theorem}[Abderrahmane and Saia--Tomazella]\label{thm-AST}
Suppose that $\{f_t\}$ is a family of isolated singularities such that the polynomial function $f_t$ is Newton non-degenerate for all small $t$. If furthermore the Milnor number $\mu(f_t)$ of $f_t$ at $\mathbf{0}$ is independent of $t$ for all small $t$, then $\{f_t\}$ is equimultiple.
\end{theorem}

Note that in this theorem, it is not required that the Newton diagram of $f_t$ (with respect to the coordinates $\mathbf{z}$) is independent of~$t$. For the definitions of Newton diagram, Newton non-degeneracy and all related topics, we refer the reader to \cite{Ko}.

In the present paper, we investigate weighted homogeneous and Newton non-degenerate singularities of hypersurfaces in the simplest case where these singularities are \emph{not} isolated, namely the case of \emph{line singularities}.
While Theorem \ref{thm-AST} easily generalizes to such a class of singularities (cf.~Theorem \ref{mt4}), extending Theorem \ref{thm-GOSh} seems to be much more challenging. Nevertheless we propose a partial generalization (cf.~Theorem \ref{mt3}).

\section{Statement of the results}\label{sect-sr}

In this section, we present our main results.

\subsection{Weighted homogeneous line singularities} 
The following theorem is a partial generalization of the theorem of Greuel and O'Shea to families $\{f_t\}$ of line singularities. As in \cite[\S4]{M7}, by this we mean that for all small $t$ the singular locus $\Sigma f_t$ of $f_t$ near the origin is given by the $z_1$-axis and the restriction of $f_t$ to the hyperplane $V(z_1)$ defined by $z_1=0$ has an isolated singularity at $\mathbf{0}$.

\begin{theorem}\label{mt3}
Suppose that $\{f_t\}$ is a family of line singularities such that the polynomial function $f_0$ is weighted homogeneous with respect to a system of positive integer weights~$(w_1,\ldots,w_n)$ satisfying the following two conditions:
\begin{enumerate}
\item[(i)]
$w_1=\mbox{\emph{min}} \{w_1,\ldots,w_n\}$;
\item[(ii)]
$w_1$ divides the weighted degree $d$ of $f_0$.
\end{enumerate} 
Also, assume that for any $t\not=0$ the polar curve $\Gamma^1_{f_t,\mathbf{z}}$ is irreducible. Under these assumptions, if furthermore the families 
\begin{equation*}
\{f_t\}
\quad\mbox{and}\quad
\{{f_t}\vert_{V(z_1)}\}
\end{equation*}
are both topologically equisingular, then they are both equimultiple.
\end{theorem}

For the definition of the polar curve $\Gamma^1_{f_t,\mathbf{z}}$, we refer the reader to \cite[Chap.~1]{M}.
Theorem \ref{mt3} is proved in Sec.~\ref{proof-mt3}.
The case where $f_0$ is homogeneous was first proved in \cite{Ey2}. In fact, Theorem 1.6 and Corollary 1.9 of \cite{Ey2} say that if $\{f_t\}$ is a topologically equisingular family of line singularities --- or even a family of line singularities with constant L\^e numbers (see \cite[Chap.~1]{M} for the definition) --- and if the polynomial function $f_0$ is homogeneous, then $\{f_t\}$ is equimultiple. In particular, this provides a complete generalization for line singularities of the theorem of  Gabri\`elov and Ku\v{s}nirenko.

\begin{remark}\label{remmt3}
By the results of Appendix \ref{ncfls}, if $n\geq 5$ then Theorem \ref{mt3} still holds true if we replace the assumption of topological equisingularity for $\{f_t\}$ and $\{{f_t}\vert_{V(z_1)}\}$ (i.e., the condition \eqref{c1p} in the list of Appendix \ref{ncfls}) by any one of the conditions \eqref{c0}, \eqref{c2}, \eqref{c3}, \eqref{c4}, \eqref{c6}, \eqref{c5} or \eqref{LG} in this list. In fact, in the proof of Theorem \ref{mt3}, we show that the condition \eqref{c1p} involved in the theorem can be replaced by the condition \eqref{c6} (and hence also by \eqref{c3}, \eqref{c4}, \eqref{c5} or \eqref{LG}) even if $n<5$.
\end{remark}

Theorem \ref{mt3} (or Remark \ref{remmt3}) may be viewed as a partial generalization for line singularities of Theorem \ref{thm-GOSh} in the sense that if $\{h_t\}$ is a family of isolated singularities in $\mathbb{C}^{n-1}$ (with coordinates $(z_2,\ldots,z_n)$) such that $\mu(h_t)=\mu(h_0)$ for all small $t$ and if $h_0$ is weighted homogeneous with respect to a system of weights~$(w_2,\ldots,w_n)$, then the corresponding family of line singularities in $\mathbb{C}^{n}$, defined by 
\begin{equation*}
f_t(z_1,z_2,\ldots,z_n):=h_t(z_2,\ldots,z_n),
\end{equation*}
satisfies the condition \eqref{c4} in the list of Appendix \ref{ncfls}. (Indeed,  under the above assumptions, for all small $t$ the polar curve $\Gamma^1_{f_t,\mathbf{z}}$ is empty --- and so the $0$th L\^e number $\lambda^0_{f_t,\mathbf{z}}$ of $f_t$ at $\mathbf{0}$ with respect to the coordinates $\mathbf{z}=(z_1,\ldots,z_n)$ is zero --- while $\mu({f_t}\vert_{V(z_1-a_1)})=\mu(h_t)$ for all small $a_1$. Now, by \cite[Sec.~1]{M7}, we know that the $1$st L\^e number $\lambda^1_{f_t,\mathbf{z}}$ coincides with the Milnor number $\mu({f_t}\vert_{V(z_1-a_1)})$ for all small $a_1\not=0$.)
Thus, since $f_0$ is weighted homogeneous with respect to the weights $(1,w_2,\ldots,w_n)$, Theorem \ref{mt3} and Remark \ref{remmt3} imply that $\{f_t\}$ (and hence $\{h_t\}$) has constant multiplicity.

\subsection{Newton non-degenerate line singularities}
Unlike the theorem of Greuel and O'Shea, the theorem of Abderrahmane, Saia and Tomazella easily extends to line singularities. More precisely we prove the following result.

\begin{theorem}\label{mt4}
Suppose that $\{f_t\}$ is a family of line singularities such that for all small $t$ the polynomial function $f_t$ is Newton non-degenerate. If furthermore this family is topologically equisingular, then it is equimultiple.
\end{theorem}

As in the case of isolated singularities, we do not require that the Newton diagram of $f_t$ is independent of~$t$. 
The theorem still holds true if we replace ``line singularities'' by ``aligned singularities,'' provided that the set of coordinates we deal with is  ``aligning''. (For the definitions of aligned singularities and aligning sets of coordinates, we refer the reader to \cite[Chap.~7]{M}.) More precisely we have the following statement.

\begin{theorem}\label{rem-mt4}
Suppose that $\{f_t\}$ is a family of $s$-dimensional aligned singularities (i.e., for all small $t$, $f_t$ has a $s$-dimensional aligned singularity at $\mathbf{0}$) such that for all small $t$ the polynomial function $f_t$ is Newton non-degenerate. Also, assume that the set of coordinates $\mathbf{z}$ is aligning for the function $f_0$ and for all the functions $f_{t_k}$, where $\{t_k\}$ is an infinite sequence such that $t_k\to 0$. Under these assumptions, if furthermore the family $\{f_t\}$ is topologically equisingular, then it is equimultiple.
\end{theorem}

As explained in the proof of \cite[Theorem 7.9]{M}, if $\{f_t\}$ is a family of aligned singularities and $\{t_k\}$ is an infinite sequence such that $t_k\to 0$, then we can use the Baire category theorem to conclude that there exists an aligning set of coordinates for  $f_0$ and for $f_{t_k}$ for all $k$. The existence of a set of coordinates which is aligning for $f_t$ for \emph{all} small $t$ is not always clear.

Actually, Theorem \ref{mt4} is a corollary of the following theorem.

\begin{theorem}\label{mt5}
Let $\{f_t\}$ be any family of equidimensional singularities (i.e., there is an integer $s$ such that the dimension  of the singular locus $\Sigma f_t$ of $f_t$ at ${\mathbf{0}}$ is equal to $s$ for all small $t$) such that the polynomial function $f_t$ is Newton non-degenerate and the L\^e numbers
\begin{equation*}
\lambda^0_{f_t,\mathbf{z}},\ldots,\lambda^s_{f_t,\mathbf{z}}
\end{equation*}
of $f_t$ at $\mathbf{0}$ with respect to the coordinates $\mathbf{z}$ are defined and independent of $t$ for all small $t$. Then the family $\{f_t\}$ is equimultiple.
\end{theorem}

Theorems \ref{mt4}, \ref{rem-mt4} and \ref{mt5} are proved in Sec.~\ref{proof-ndc}

\section{Examples}\label{appli}

Theorems \ref{mt3}--\ref{mt5} may be quite useful to decide whether certain families of hypersurfaces with line singularities are \emph{not} topologically  equisingular --- a question which is, in general, extremely difficult to answer. 

For example, consider the family defined by 
\begin{equation*}
f_t(z_1,z_2,z_3)=z_1^4z_2^2+z_2^4+z_3^4+tz_1z_2^2+t^2z_1^2z_2^2.
\end{equation*}
A priori, it is far from being obvious to decide whether this family is topologically  equisingular or not. However this easily follows from Theorem \ref{mt3}. Indeed, the polynomial function $f_0(z_1,z_2,z_3)=z_1^4z_2^2+z_2^4+z_3^4$ is weighted homogeneous with respect to the weights $(w_1,w_2,w_3)=(1,2,2)$, the singular locus $\Sigma f_t$ of $f_t$ near $\mathbf{0}$ is given by the $z_1$-axis, and the restriction $f_t\vert_{V(z_1)}$ has an isolated singularity at~$\mathbf{0}$ with a Milnor number independent of $t$. In particular, since $f_t\vert_{V(z_1)}$ is a function of two variables, the family $\{f_t\vert_{V(z_1)}\}$ is topologically equisingular. Finally, an easy computation shows that 
\begin{align*}
\Gamma^1_{f_t,\mathbf{z}} 
= V(2z_1^4+4z_2^2+2tz_1+2t^2z_1^2,z_3^3),
\end{align*}
which is clearly irreducible for all small $t\not=0$. 
Since the family $\{f_t\}$ is not equimultiple, it follows from Theorem \ref{mt3} that it is not topologically  equisingular.

\begin{remark}
By \cite[Corollary 3.7]{ER1}, we know that if $\{f_t\}$ is a non-equimultiple family of line singularities of the form $f_t(\mathbf{z})=f_0(\mathbf{z})+\xi(t)g(\mathbf{z})$, where $\xi\colon (\mathbb{C},0)\to (\mathbb{C},0)$ is a non-constant polynomial function and $g\colon (\mathbb{C}^n,\mathbf{0})\to (\mathbb{C},0)$ is any polynomial function, then $\{f_t\}$ is not topologically  equisingular. The above example is not a consequence of this result.
\end{remark}

More generally, suppose that $\{f_t\}$ is a family of line singularities such that $f_0$ is weighted homogeneous with respect to  weights~$(w_1,\ldots,w_n)$ satisfying the conditions (i) and (ii) of Theorem \ref{mt3}. Also, assume that for all small $t\not=0$, the polar curve $\Gamma^1_{f_t,\mathbf{z}}$ is irreducible and any monomial of $f$ that contains a non-zero power of $t$ also contains a non-zero power of $z_1$. Under these assumptions, if the multiplicity of $f_t$ jumps at $t=0$, then the family $\{f_t\}$ is not topologically equisingular.

Let us now give an example in the Newton non-degenerate case. For instance, we easily check that the family $\{f_t\}$ defined by 
\begin{equation*}
f_t(z_1,z_2,z_3)=z_1^2z_2^2+z_2^4+z_2^5+z_3^5+tz_1z_2^2+t^2z_1^2z_2^3
\end{equation*}
is a family of line singularities such that for all $t$ small enough $f_t$ is  Newton non-degenerate. (Note that the Newton diagram of $f_t$ changes at $t=0$.) Since $\{f_t\}$ is not equimultiple, Theorem \ref{mt4} implies that it is not topologically equisingular.

\section{Proof of Theorem \ref{mt3}}\label{proof-mt3}

Roughly, the idea of the proof is to use the Iomdine-L\^e-Massey formula (cf.~\cite[Theorem 4.5]{M}) to reduce the problem to the case of isolated singularities and then apply the Greuel-O'Shea theorem (cf.~Theorem \ref{thm-GOSh}). In our case, the Iomdine-L\^e-Massey formula says that for all but a finite number of non-zero complex numbers $a(t)$, the function $f_t+a(t)z_1^{\rho_t}$ has an isolated singularity at $\mathbf{0}$ with Milnor number equal to  $\lambda^0_{f_t,\mathbf{z}}+(\rho_t-1)\lambda^1_{f_t,\mathbf{z}}$,
where $\rho_t$ is the maximum value between the number $2$ and the maximum ``polar ratio'' for $f_t$ at~$\mathbf{0}$, and where $\lambda^0_{f_t,\mathbf{z}}$, $\lambda^1_{f_t,\mathbf{z}}$ are the L\^e numbers of $f_t$ at $\mathbf{0}$ with respect to the coordinates $\mathbf{z}$. 
The polar ratios, which are essential in the proof, are described in Lemmas \ref{lemma-int} and \ref{lemma1} below.
 In particular, note that the weighted homogeneity of the function $f_0+a(0)z_1^{\rho_0}$ is controlled by the maximum polar ratio of $f_0$.

Now let us go into the details. First, observe that if $d/w_1=1$, then the polynomial function $f_0$ is homogeneous, and the result follows from \cite[Theorem 1.6 and Corollary 1.9]{Ey2}. 
From now on, we assume that $d/w_1\geq 2$. By the results of Appendix \ref{ncfls}, the L\^e numbers $\lambda^1_{f_t,\mathbf{z}}$, $\lambda^0_{f_t,\mathbf{z}}$ and the polar number $\gamma^1_{f_t,\mathbf{z}}$ are defined and independent of $t$ for all small~$t$ (cf.~Appendix \ref{a6}).\footnote{Note that for line singularities, the only possible non-zero L\^e numbers are precisely $\lambda^0_{f_t,\mathbf{z}}$ and $\lambda^1_{f_t,\mathbf{z}}$; all the other L\^e numbers $\lambda^k_{f_t,\mathbf{z}}$ for $2\leq k\leq n-1$ are defined and equal to zero (cf.~\cite{M}).} In particular, the polar curve $\Gamma^1_{f_t,\mathbf{z}}$ is purely $1$-dimensional or empty at~$\mathbf{0}$ and the intersection $\Gamma^1_{f_t,\mathbf{z}}\cap V(z_1)$ is  $0$-dimensional or empty at $\mathbf{0}$. 
Let $\eta$ be an irreducible component of $\Gamma^1_{f_t,\mathbf{z}}$ (with its reduced structure) which is not empty at $\mathbf{0}$. 
By \cite[Definition 4.1]{M}, the \emph{polar ratio} $\rho(\eta)$ of $\eta$ at $\mathbf{0}$ is defined by
\begin{equation*}
\rho(\eta):=\frac{([\eta]\cdot [V(f_t)])_\mathbf{0}}
{([\eta]\cdot [V(z_1)])_\mathbf{0}}.
\end{equation*}
Here, $([\eta]\cdot [V(f_t)])_\mathbf{0}$ denotes the intersection number at $\mathbf{0}$ of the analytic cycles $[\eta]$ and $[V(f_t)]$ associated to the schemes $\eta$ and $V(f_t)$ respectively. Similarly for $([\eta]\cdot [V(z_1)])_\mathbf{0}$. A polar ratio for $f_t$ at $\mathbf{0}$ is any one of the polar ratios of any component $\eta$ of $\Gamma^1_{f_t,\mathbf{z}}$ such that $\eta$ is not empty at $\mathbf{0}$. (If $\Gamma^1_{f_t,\mathbf{z}}$ is empty at $\mathbf{0}$, then we say that ``the maximum polar ratio for $f_t$ at $\mathbf{0}$ is $1$.'') 

The polar ratios for $f_0$ and $f_t$ ($t\not=0$) at $\mathbf{0}$ are described in Lemmas \ref{lemma-int} and \ref{lemma1} respectively. Before to state these lemmas, note that we may assume that $\gcd(w_1,\ldots,w_n)=1$. (Indeed, if $\gcd(w_1,\ldots,w_n)=w>1$, we take $w_i':=w_i/w$ ($1\leq i\leq n$). Then, clearly, $\gcd(w'_1,\ldots,w'_n)=1$, $w'_1=\mbox{min}\{w'_1,\ldots,w'_n\}$, and $f_0$ is weighted homogeneous with respect to $(w'_1,\ldots,w'_n)$.)

\begin{lemma}[see also Proposition 3.10 of \cite{MS}]\label{lemma-int}
If $\eta$ is an irreducible component of $\Gamma^1_{f_0,\mathbf{z}}$ which is not empty at $\mathbf{0}$, then
\begin{equation*}
([\eta]\cdot [V(f_0)])_\mathbf{0}=d
\quad\mbox{and}\quad
([\eta]\cdot [V(z_1)])_\mathbf{0}=w_1.
\end{equation*}
In particular, if $\gamma^1_{f_0,\mathbf{z}}\not=0$ (equivalently, if $\Gamma^1_{f_0,\mathbf{z}}$ is not empty at $\mathbf{0}$), then the polar ratios for $f_0$ at $\mathbf{0}$ are equal to $d/w_1$.
\end{lemma}

\begin{proof}
First, observe that in any neighbourhood of the origin, $\eta$ has a point of the form $(a_1,\ldots,a_n)$ with $a_1\not=0$. (Otherwise, there is a neighbourhood  in which all the points of $\eta$ are of the form $(0,a_2,\ldots,a_n)$. Thus, in such a neighbourhood, $\eta\subseteq V(z_1)$ and $\dim_{\mathbf{0}}(\eta\cap V(z_1))=\dim_{\mathbf{0}}\eta=1$ --- a contradiction.) 

We claim that $\eta\not\subseteq V(f_0)$.
Indeed, since $f_0$ is weighted homogeneous with respect to $(w_1,\ldots,w_n)$ and since $\gcd(w_1,\ldots,w_n)=1$, we may pick a parametrization of $\eta$ of the form 
\begin{align*}
s\mapsto \phi(s):=(a_1 s^{w_1},\ldots, a_n s^{w_n}).
\end{align*}
Clearly, for any $i\geq 2$,
\begin{equation*}
\frac{\partial f_0}{\partial z_i}(\phi(s)) =
\frac{\partial f_0}{\partial z_i}(a_1 s^{w_1},\ldots, a_n s^{w_n}) =
s^{d-w_i} \frac{\partial f_0}{\partial z_i}(a_1,\ldots, a_n) = 0.
\end{equation*}
 Thus, if $f_0\circ\phi$ identically vanishes, then, for all $s$,
\begin{equation*}
0=(f_0\circ\phi)'(s)=a_1w_1s^{w_1-1} \frac{\partial f_0}{\partial z_1}(\phi(s)),
\end{equation*}
and hence $\frac{\partial f_0}{\partial z_1}(\phi(s))=0$. 
It follows that $\eta$ is contained in $\Sigma f_0$ --- a contradiction.

Now, since $\eta\not\subseteq V(f_0)$, a classical result in intersection theory (cf.~\cite{Fulton} or \cite[Chap.~1]{M}) shows that 
\begin{align*}
([\eta]\cdot [V(f_0)])_\mathbf{0} = \mbox{ord}_{0} (f_0\circ \phi(s))=d,
\end{align*} 
where $\mbox{ord}_{0} (f_0\circ \phi(s))$ is the order of $f_0\circ \phi(s)$ at $0$.

A similar argument shows $([\eta]\cdot [V(z_1)])_\mathbf{0}=w_1$.
\end{proof}

\begin{lemma}\label{lemma1}
For any sufficiently small $t\not=0$, the unique polar ratio $\rho(\Gamma^1_{f_t,\mathbf{z}})$ for $f_t$ at $\mathbf{0}$ is equal to $d/w_1$ if $\gamma^1_{f_t,\mathbf{z}}\not=0$.
\end{lemma}

\begin{proof}
By \cite[Proposition 1.23]{M}, $([\Gamma^1_{f_t,\mathbf{z}}]\cdot [V(f_t)])_\mathbf{0}= \gamma^1_{f_t,\mathbf{z}}+\lambda^0_{f_t,\mathbf{z}}$.
Thus, writing the cycle $[\Gamma^1_{f_0,\mathbf{z}}]$ as $\sum_\eta k_\eta [\eta]$, where the sum is taken over all the components $\eta$ of $\Gamma^1_{f_0,\mathbf{z}}$, it follows from Lemma \ref{lemma-int} that for all small $t\not=0$:
\begin{align*}
\rho(\Gamma^1_{f_t,\mathbf{z}}) & =\frac{\gamma^1_{f_t,\mathbf{z}}+\lambda^0_{f_t,\mathbf{z}}}{\gamma^1_{f_t,\mathbf{z}}}=
\frac{\gamma^1_{f_0,\mathbf{z}}+\lambda^0_{f_0,\mathbf{z}}}{\gamma^1_{f_0,\mathbf{z}}}
=\frac{([\Gamma^1_{f_0,\mathbf{z}}]\cdot [V(f_0)])_\mathbf{0}}{([\Gamma^1_{f_0,\mathbf{z}}]\cdot [V(z_1)])_\mathbf{0}}\\
& =\frac{\sum k_\eta([\eta]\cdot [V(f_0)])_\mathbf{0}}
{\sum k_\eta([\eta]\cdot [V(z_1)])_\mathbf{0}}=\frac{d}{w_1}.\qedhere
\end{align*}
\end{proof}

Now we precisely know the value of the polar ratios, we can use the Iomdine-L\^e-Massey formula to show the following lemma.

\begin{lemma}\label{lemma43}
If $\gamma^1_{f_t,\mathbf{z}}\not=0$, then there exists a countable subset $E\subseteq\mathbb{C}$ such that for any non-zero element $a\notin E$ and any sufficiently small $t$, the function
\begin{equation*}
f_t+a z_1^{d/w_1}
\end{equation*}
has an isolated singularity at $\mathbf{0}$ with a Milnor number independent of $t$.
\end{lemma}

\begin{proof}
As in the proof of Theorem 7.9 of \cite{M}, we first show the desired property for an infinite sequence $\{t_k\}$ of the parameter $t$ approaching $0$ and then we deduce the property for all small $t$ using the upper-semicontinuity of the Milnor number.

So, let $\{t_k\}$ be an infinite sequence such that $t_k\to 0$. In the light of Lemmas 4.1 and 4.2 and since $d/w_1\geq 2$, by applying the Iomdine-L\^e-Massey formula (cf.~\cite[Theorem 4.5]{M}) to the function $f_{t_k}$ for a given fixed $t_k$, we obtain that for all but a finite number of non-zero complex numbers $a(t_k)$, depending on that $t_k$, the function $f_{t_k}+a(t_k) z_1^{d/w_1}$ has an isolated singularity at~$\mathbf{0}$ with Milnor number 
\begin{equation}\label{Mn}
\mu(f_{t_k}+a(t_k) z_1^{d/w_1})=\lambda^0_{f_{t_k},\mathbf{z}}+\biggl(\frac{d}{w_1}-1\biggr)\lambda^1_{f_{t_k},\mathbf{z}}.
\end{equation}
The set $E(t_k)$ of the excluded values of $a(t_k)$ consists of those non-zero complex numbers $a(t_k)$ which make the lowest degree terms of
\begin{equation*}
\biggl(\frac{\partial f_{t_k}}{\partial z_1}\biggr)\big\vert_{\phi_\eta(s)}
\quad\mbox{and}\quad
\biggl(\frac{d}{w_1} a(t_k) z_1^{\frac{d}{w_1}-1}\biggr)\big\vert_{\phi_\eta(s)}
\end{equation*}
add up to zero, where $\eta$ is an irreducible component of $\Gamma^1_{f_{t_k},\mathbf{z}}$ and $\phi_\eta(s)$ is a param\-etri\-zation of $\eta$ (see the proof of \cite[Lemma 4.3]{M}). 

A similar property also holds true if we replace $t_k$ by $t=0$. If $E(0)$ denotes the corresponding set of excluded values, then for any non-zero element $a$ which is not contained in the countable set $E:=E(0)\cup\bigcup_k E(t_k)$, the functions 
\begin{equation*}
f_0+a z_1^{d/w_1}
\quad\mbox{and}\quad
f_{t_k}+a z_1^{d/w_1}
\end{equation*}
(all $k$) have an isolated singularity at $\mathbf{0}$ and the same Milnor number \eqref{Mn} (remind that the L\^e numbers are independent of $t$).
This, together with the upper-semicontinuity of the Milnor number, imply that for all $t$ small enough the function $f_{t}+a z_1^{d/w_1}$ has an isolated singularity at $\mathbf{0}$ and the same Milnor number as $f_0+a z_1^{d/w_1}$.
\end{proof}

Since for $t=0$ the function $f_0+a z_1^{d/w_1}$  is weighted homogeneous, the Greuel-O'Shea theorem (cf.~Theorem \ref{thm-GOSh}) says that 
\begin{equation}\label{mainequality1}
\mbox{ord}_{\mathbf{0}}(f_t+az_{1}^{d/w_{1}}) 
= \mbox{ord}_{\mathbf{0}}(f_0+az_{1}^{d/w_{1}}),
\end{equation}
where $\mbox{ord}_{\mathbf{0}}(f_t+az_{1}^{d/w_{1}})$ is the order of $f_t+az_{1}^{d/w_{1}}$ at $\mathbf{0}$.

If $\mbox{ord}_{\mathbf{0}}(f_0)\leq d/w_{1}$, then (\ref{mainequality1}) immediately implies $\mbox{ord}_{\mathbf{0}}(f_t)=\mbox{ord}_{\mathbf{0}}(f_0)$ and the theorem is proved.
Now we claim that we always have $\mbox{ord}_{\mathbf{0}}(f_0)\leq d/w_{1}$. 
Indeed, take any monomial $c\, z_1^{\alpha_1}\cdots z_n^{\alpha_n}$ ($c$ constant) of the initial polynomial $\mbox{in} (f_0)$ of $f_0$. Then we have 
$\sum_{1\leq i\leq n} \alpha_i=\deg \mbox{in} (f_0) = \mbox{ord}_{\mathbf{0}}(f_0)$, and since $f_0$ (and hence $\mbox{in} (f_0)$) is weighted homogeneous with respect to the weights $(w_1,\ldots,w_n)$, we also have $\sum_{1\leq i\leq n}\alpha_i w_i=d$.
As $w_{1}$ is the smallest weight, it follows that
\begin{equation*}
d=\sum_{1\leq i\leq n}\alpha_i w_i \geq \sum_{1\leq i\leq n}\alpha_i w_{1} = w_{1} \sum_{1\leq i\leq n}\alpha_i = w_{1} \mbox{ord}_{\mathbf{0}}(f_0).
\end{equation*}
This completes the proof of Theorem \ref{mt3} when  $\gamma^1_{f_t,\mathbf{z}}\not=0$.

\begin{remark}
That $\{{f_t}\vert_{V(z_1)}\}$ is equimultiple immediately follows from the Greuel-O'Shea theorem. Indeed, by the assumption, $\{{f_t}\vert_{V(z_1)}\}$ is a topologically equisingular family of isolated singularities with ${f_0}\vert_{V(z_1)}$ weighted homogeneous.
\end{remark}

If $\gamma^1_{f_t,\mathbf{z}}=0$, then $\mathbf{0}\not\in\Gamma^1_{f_t,\mathbf{z}}$ (i.e., $\Gamma^1_{f_t,\mathbf{z}}$ is empty at $\mathbf{0}$), and the maximum polar ratio for $f_t$ at $\mathbf{0}$ is $1$. Since $d/ w_1> 1$, the Iomdine-L\^e-Massey formula shows that for all small $t$ and all non-zero complex numbers $a$, the function $f_t+az_{1}^{d/w_{1}}$ has an isolated singularity at~$\mathbf{0}$ with a Milnor number independent of $t$. Then we conclude exactly as above using the Greuel-O'Shea theorem.

\section{Proofs of Theorems \ref{mt4}--\ref{mt5}}\label{proof-ndc}

Again, the idea is to use the Iomdine-L\^e-Massey formula to reduce the problem to the case of isolated singularities and then apply the theorem of Abderrahmane, Saia and Tomazella (cf.~Theorem \ref{thm-AST}). This time, in order to control the Newton non-degeneracy, we use a result of Brzostowski and Oleksik \cite[Lemma 3.7]{BO}. This lemma says that if $h\colon (\mathbb{C}^n,\mathbf{0})\to(\mathbb{C},0)$ is a Newton non-degenerate holomorphic function with a (possibly non-isolated) singularity at $\mathbf{0}$, then there exists a (non unique) constant $C(h)$ such that for any $j>C(h)$, the function $h+z_i^{j}$ is Newton non-degenerate too. The constant $C(h)$ constructed in the proof of the lemma only depends on the Newton diagram of $h$. In particular, in our case, $C(f_t)$ does not depend on $t$ for all $t\not=0$ but we may have $C(f_t)\not=C(f_0)$, so that Newton non-degeneracy is guaranteed if $j>\max\{C(f_t),C(f_0)\}$.

\enlargethispage{1mm}
The proofs below are inspired by the proof of Theorem 7.9 in \cite{M}.

\subsection{Proof of Theorem \ref{mt5}}
By the Brzostowski-Oleksik lemma \cite[Lemma 3.7]{BO} and the Iomdine-L\^e-Massey formula \cite[Theorem 4.5 and Corollary~4.6]{M}, for all integers $0\ll j_1\ll j_2\ll \ldots\ll j_{s}$, the function
\begin{equation*}
f_t+z_1^{j_1}+\ldots +z_{s}^{j_{s}}
\end{equation*}
is Newton non-degenerate and has an isolated singularity at~$\mathbf{0}$ with a Milnor number independent of $t$ provided that $t$ is small enough. Then Theorem \ref{thm-AST} implies that the family $\{f_t+z_1^{j_1}+\ldots +z_{s}^{j_{s}}\}$ is equimultiple. Since $j_1,\ldots,j_s$ are arbitrarily large, we deduce that $\{f_t\}$ is equimultiple too.

\subsection{Proof of Theorem \ref{mt4}}
Theorem \ref{mt4} is a corollary of Theorem \ref{mt5}.
Indeed, by \cite[Remark 1.29]{M}, the partition of $V(f_t)$ given by
$\mathscr{S}_t:=\bigl\{V(f_t)\setminus\Sigma f_t,\Sigma f_t\setminus\{\mathbf{0}\}, \{\mathbf{0}\}\bigr\}$ is a \emph{good stratification} for $f_t$ in a neighbourhood of $\mathbf{0}$, and the hyperplane $V(z_1)$ is a prepolar slice for $f_t$ at $\mathbf{0}$ with respect to $\mathscr{S}_t$ for all small $t$. 
In particular, combined with \cite[Proposition~1.23]{M}, this implies that the L\^e numbers 
\begin{equation*}
\lambda^0_{f_t,\mathbf{z}}
\quad\mbox{and}\quad
\lambda^1_{f_t,\mathbf{z}}
\end{equation*}
 of $f_t$ at~$\mathbf{0}$ with respect to the coordinates $\mathbf{z}$ are defined. 
Now, from \cite[Sec.~4]{M7}, we know that $\lambda^0_{f_t,\mathbf{z}}$ and $\lambda^1_{f_t,\mathbf{z}}$ are independent of $t$ for all small $t$ if and only if the generic Milnor number $\mathring{\mu}(f_t)$ (i.e., the Milnor number $\mu({f_t}\vert_{V(z_1-a_1)})$, which is independent of $a_1$ for all small $a_1\not=0$) and the reduced Euler characteristic of the Milnor fibre of $f_t$ at~$\mathbf{0}$ are both independent of $t$ for all small $t$. We also know that for line singularities, $\mathring{\mu}(f_t)$ is an invariant of the ambient topological type of $V(f_t)$ at~$\mathbf{0}$ (cf.~\cite[Sec.~1]{M7}).
Since the reduced Euler characteristic is a topological invariant too, and since our family $\{f_t\}$ is topologically equisingular, it follows that the L\^e numbers $\lambda^0_{f_t,\mathbf{z}}$ and $\lambda^1_{f_t,\mathbf{z}}$ are independent of $t$ for all small $t$. Then Theorem \ref{mt5} implies that $\{f_t\}$ is equimultiple.

\subsection{Proof of Theorem \ref{rem-mt4}}
Since the set of coordinates $\mathbf{z}$ is aligning for $f_0$ and $f_{t_k}$,  the L\^e numbers $\lambda^q_{f_{t_k},\mathbf{z}}$ and $\lambda^q_{f_{0},\mathbf{z}}$  are defined and \cite[Corol\-lary~7.8]{M} implies
\begin{equation*}
\lambda^q_{f_{t_k},\mathbf{z}} = \lambda^q_{f_{0},\mathbf{z}}
\end{equation*}
 for all $0\leq q\leq s$ and all $k$ sufficiently large. Thus, by \cite[Lemma 3.7]{BO} and \cite[Theorem 4.5 and Corollary~4.6]{M} again, for all integers $0\ll j_1\ll j_2\ll \ldots\ll j_{s}$, the functions 
\begin{equation*}
f_0+z_1^{j_1}+\ldots +z_{s}^{j_{s}}
\quad\mbox{and}\quad
f_{t_k}+z_1^{j_1}+\ldots+z_{s}^{j_{s}}
\end{equation*}
are Newton non-degenerate and have an isolated singularity at~$\mathbf{0}$ with the same Milnor number, provided that $k$ is large enough. By the upper-semicontinuity of the Milnor number, this implies that for all small $t$, the function $f_t+z_1^{j_1}+\ldots +z_{s}^{j_{s}}$ has an isolated singularity at $\mathbf{0}$ with the same Milnor number as $f_0+z_1^{j_1}+\ldots +z_{s}^{j_{s}}$. Clearly, it also implies that $f_t+z_1^{j_1}+\ldots +z_{s}^{j_{s}}$ is Newton non-degenerate for all small $t$. Then we conclude as in the proof of Theorem \ref{mt5} using Theorem \ref{thm-AST}.

\appendix

\section{Numerical control of the topological type in a family of line singularities}\label{ncfls}

To understand better the assumptions used in Theorem \ref{mt3} and Remark \ref{remmt3}, we discuss here how the topological type in a family of line singularities can be controlled by numerical invariants.
All the results presented in this appendix are well known (or are immediate consequences of well known theorems). 

Suppose that $\{f_t\}$ is a family of line singularities. 
We recall that in this case the L\^e numbers $\lambda^0_{f_t,\mathbf{z}}$, $\lambda^1_{f_t,\mathbf{z}}$ and the polar number $\gamma^1_{f_t,\mathbf{z}}$ of $f_t$ at~$\mathbf{0}$ with respect to the coordinates $\mathbf{z}$ are defined (cf.~\cite{M}).
Then we want to discuss the relationships between the following conditions.

\begin{enumerate}
\item\label{c0}
The singular locus $\Sigma f$ of the underlying function $f$ defining the family $\{f_t\}$ is given by $V(z_2,\ldots,z_n):=\{(t,\mathbf{z})\in\mathbb{C}\times\mathbb{C}^n\, ;\, z_i=0 \mbox{ for } i\geq 2\}$ and the Milnor number $\mu(f\vert_{V(t-t_0,z_1-a_1)})$ is independent of the point $(t_0,a_1,\mathbf{0})$ in some open neighbourhood of $(0,0,\mathbf{0})$.
\item\label{c2}
For any sufficiently small $a_1$, the families $\{f_t\}$ and $\{{f_t}\vert_{V(z_1-a_1)}\}$ are topologically  equisingular.
\item\label{c1p}
The families $\{f_t\}$ and $\{{f_t}\vert_{V(z_1)}\}$ are topologically  equisingular.
\item\label{c1}
The family $\{f_t\}$ is topologically  equisingular.
\item\label{c3}
The generic Milnor number $\mathring\mu({f_t})$ and the Milnor numbers $\mu({f_t}\vert_{V(z_1)})$ and $\mu({f_t}+z_1^j)$ are constant for any sufficiently large $j$.
\item\label{c3p}
The generic Milnor number $\mathring\mu({f_t})$ and the Milnor number $\mu({f_t}+z_1^j)$ are constant for any sufficiently large $j$.
\item\label{c4}
The Milnor number $\mu({f_t}\vert_{V(z_1)})$ and the L\^e numbers $\lambda^0_{f_t,\mathbf{z}}$ and $\lambda^1_{f_t,\mathbf{z}}$ are constant.
\item\label{c6}
The polar number $\gamma^1_{f_t,\mathbf{z}}$ and the L\^e numbers $\lambda^0_{f_t,\mathbf{z}}$ and $\lambda^1_{f_t,\mathbf{z}}$ are constant.
\item\label{c5}
The sum $\gamma^1_{f_t,\mathbf{z}}+\lambda^0_{f_t,\mathbf{z}}$ and the L\^e number $\lambda^1_{f_t,\mathbf{z}}$ are constant.
\item\label{c7}
The L\^e numbers $\lambda^0_{f_t,\mathbf{z}}$ and $\lambda^1_{f_t,\mathbf{z}}$ are constant.
\item\label{c7bis}
The L\^e number $\lambda^0_{f_t,\mathbf{z}}$ and the generic Milnor number $\mathring\mu({f_t})$ are constant.
\item\label{c8}
The generic Milnor number $\mathring\mu({f_t})$ and the reduced Euler characteristic $\widetilde\chi(F_{f_t,\mathbf{0}})$ of the Milnor fibre $F_{f_t,\mathbf{0}}$ of $f_t$ at $\mathbf{0}$ are constant.
\item\label{LG}
The generic Milnor number $\mathring\mu({f_t})$ and the number
\begin{equation*}
m(f_t+z_1^j):=
\dim_{\mathbb{C}}\mathscr{O}_{n,\mathbf{0}}\big/(f_t+z_1^j,J(f_t+z_1^j,z_1))
\end{equation*}
are constant. (Here, $J(f_t+z_1^j,z_1)$ denotes the ideal generated by the determinants of all $2$-minors of the corresponding Jacobian matrix, that is, the ideal generated by $\frac{\partial (f_t+z_1^j)}{\partial z_i}$ for all $2\leq i\leq n$.)
\item\label{cr1}
The L\^e numbers $\lambda^0_{f_t,\mathbf{z}}$ and $\lambda^1_{f_t,\mathbf{z}}$ and the polar number $\gamma^1_{f_t,\mathbf{z}}$ and $\gamma^2_{f_t,\mathbf{z}}$ are constant.\footnote{For the definition of $\gamma^2_{f_t,\mathbf{z}}$, which we have not encountered yet, we also refer the reader to \cite{M}.}
\item\label{cmf}
The diffeomorphism type of the Milnor fibration of $f_t$ at $\mathbf{0}$ is constant.
\item\label{c00}
The family $\{f_t\}$ is Whitney equisingular (i.e., there is a \emph{Whitney stratification} of $V(f)$ in a neighbourhood of $\mathbf{0}$ such that the $t$-axis $\mathbb{C}\times\{\mathbf{0}\}$ is a stratum).
\item\label{equimult}
The family $\{f_t\}$ is equimultiple.
\end{enumerate}

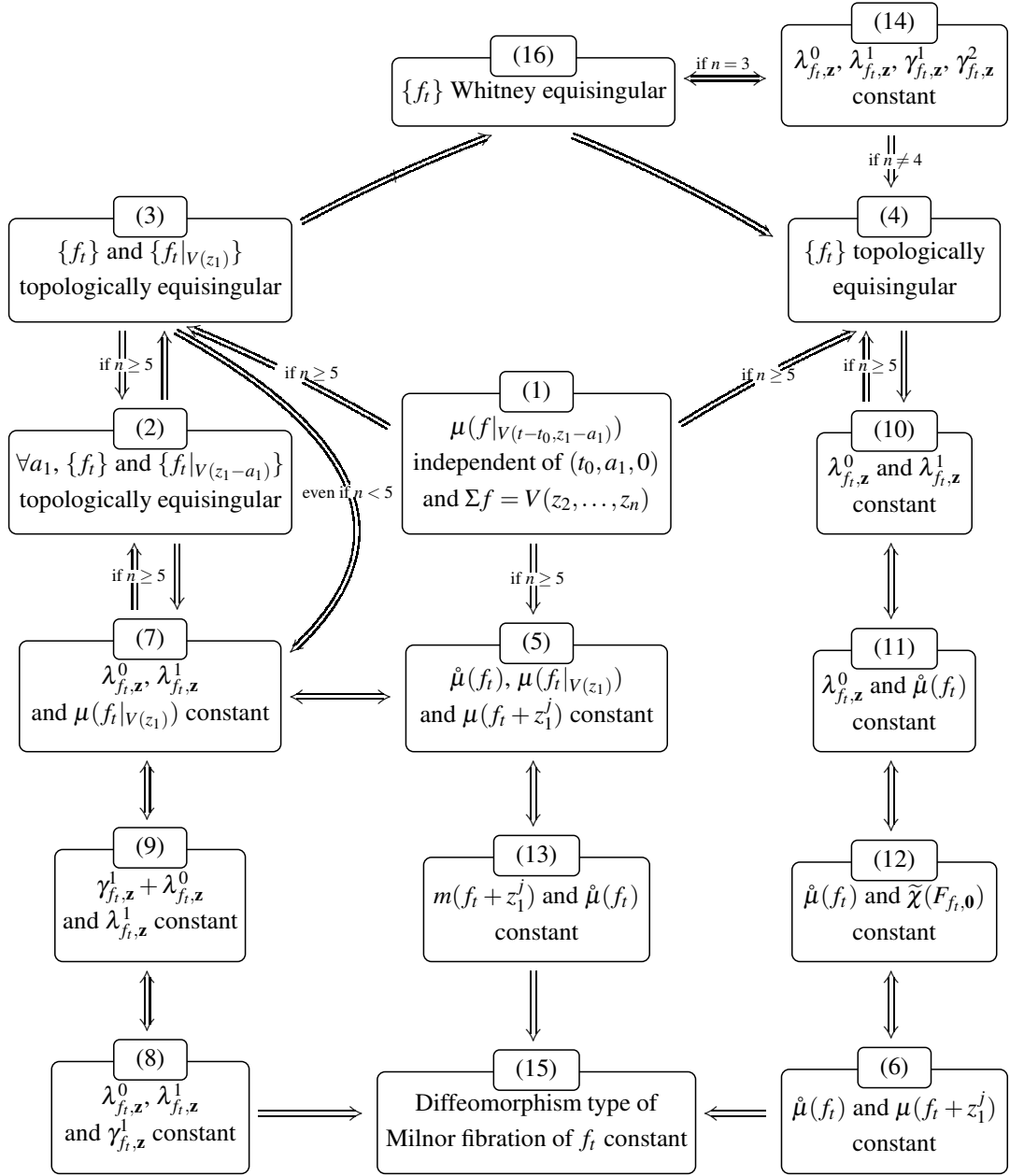
\begin{figure}[t]
\centering
\onehalfspacing
\small{
\begin{eqnarray*}
\xymatrix{
& 
\mbox{\begin{mybox}[colback=white,width = 4cm]{\textbf{\mbox{\eqref{c00}}}}
\centerline{$\{f_t\}$ Whitney equisingular}
\end{mybox}}
\ar@/_0.5pc/@{<=}[dl] |-{\SelectTips{cm}{}\object@{/}}
\ar@/^0.5pc/@{=>}[dr] 
\ar@{<=>}@<17pt>[r]^{\mbox{\tiny $\quad$ if $n=3$}}
& 
\mbox{
\begin{mybox}[colback=white,width = 3.2cm]{\textbf{\mbox{\eqref{cr1}}}}
\centerline{$\lambda^0_{f_t,\mathbf{z}}$, $\lambda^1_{f_t,\mathbf{z}}$, $\gamma^1_{f_t,\mathbf{z}}$, $\gamma^2_{f_t,\mathbf{z}}$} 
\centerline{constant}
\end{mybox}}
\\ 
\mbox{
\begin{mybox}[colback=white,width = 4cm]{\textbf{\mbox{\eqref{c1p}}}}
\centerline{$\{f_t\}$ and $\{{f_t}\vert_{V(z_1)}\}$}\centerline{topologically equisingular}
\end{mybox}}
\ar@/^7pc/@{=>}[dd]|-{\mbox{\tiny \hskip -0.5mm even if $n<5$}}
& 
& \mbox{
\begin{mybox}[colback=white,width = 3cm]{\textbf{\mbox{\eqref{c1}}}}
\centerline{$\{f_t\}$ topologically}
\centerline{equisingular}
\end{mybox}}
\ar@{<=}[u]|-{\mbox{\tiny \ \ if $n\not=4$}}
\\
\mbox{
\begin{mybox}[colback=white,width = 4cm]{\textbf{\mbox{\eqref{c2}}}}
\centerline{$\forall a_1$, $\{f_t\}$ and $\{{f_t}\vert_{V(z_1-a_1)}\}$}\centerline{topologically equisingular}
\end{mybox}}
\ar@{=>}@<-7pt>[u] 
\ar@{<=}@<10pt>[u]|-{\mbox{\tiny \ \ if $n\geq 5$}}
& 
\mbox{
\begin{mybox}[colback=white,width = 3.8cm]{\textbf{\mbox{\eqref{c0}}}}
\centerline{$\mu(f\vert_{V(t-t_0,z_1-a_1)})$}
\centerline{independent of $(t_0,a_1,0)$}
\centerline{and $\Sigma f=V(z_2,\ldots,z_n)$}
\end{mybox}}
\ar@/^0.5pc/@{=>}[ur]|-{\mbox{\tiny \ if $n\geq 5$}}
\ar@/_0.5pc/@{=>}@<-2pt>[ul]|-{\mbox{\tiny $\qquad$ \ if $n\geq 5$}}
\ar@{=>}[d]|-{\mbox{\tiny \ \ if $n\geq 5$}}
& 
\mbox{
\begin{mybox}[colback=white,width = 2.2cm]{\textbf{\mbox{\eqref{c7}}}}
\centerline{$\lambda^0_{f_t,\mathbf{z}}$ and $\lambda^1_{f_t,\mathbf{z}}$} \centerline{constant} 
\end{mybox}}
\ar@{<=}@<-5pt>[u] 
\ar@{=>}@<10pt>[u]|-{\mbox{\tiny \ \ if $n\geq 5$}}
\\ 
\mbox{
\begin{mybox}[colback=white,width = 3.7cm]{\textbf{\mbox{\eqref{c4}}}}
\centerline{$\lambda^0_{f_t,\mathbf{z}}$, $\lambda^1_{f_t,\mathbf{z}}$} 
\centerline{and $\mu({f_t}\vert_{V(z_1)})$ constant}
\end{mybox}}
\ar@{<=>}@<20pt>[r]
\ar@{<=>}[d]
\ar@{<=}@<-12pt>[u] 
\ar@{=>}@<5pt>[u]|-{\mbox{\tiny \ \ if $n\geq 5$}}
& 
\mbox{
\begin{mybox}[colback=white,width = 3.7cm]{\textbf{\mbox{\eqref{c3}}}}
\centerline{$\mathring\mu({f_t})$, $\mu({f_t}\vert_{V(z_1)})$}\centerline{and $\mu({f_t}+z_1^j)$ constant} 
\end{mybox}}
& 
\mbox{
\begin{mybox}[colback=white,width = 2.3cm]{\textbf{\mbox{\eqref{c7bis}}}}
\centerline{$\lambda^0_{f_t,\mathbf{z}}$ and $\mathring\mu({f_t})$} \centerline{constant}
\end{mybox}}
\ar@{<=>}[u]
\\ 
\mbox{
\begin{mybox}[colback=white,width = 2.7cm]{\textbf{\mbox{\eqref{c5}}}}
\centerline{$\gamma^1_{f_t,\mathbf{z}}+\lambda^0_{f_t,\mathbf{z}}$}\centerline{and $\lambda^1_{f_t,\mathbf{z}}$ constant}
\end{mybox}}
\ar@{<=>}[d]
& 
\mbox{
\begin{mybox}[colback=white,width = 3.2cm]{\textbf{\mbox{\eqref{LG}}}}
\centerline{$m({f_t}+z_1^j)$ and $\mathring\mu({f_t})$} \centerline{constant}
\end{mybox}}
\ar@{<=>}[u]
& 
\mbox{
\begin{mybox}[colback=white,width = 2.9cm]{\textbf{\mbox{\eqref{c8}}}}
\centerline{$\mathring\mu({f_t})$ and $\widetilde\chi(F_{f_t,\mathbf{0}})$}\centerline{constant}
\end{mybox}}
\ar@{<=>}[u]
\\ 
\mbox{
\begin{mybox}[colback=white,width = 2.8cm]{\textbf{\mbox{\eqref{c6}}}}
\centerline{$\lambda^0_{f_t,\mathbf{z}}$, $\lambda^1_{f_t,\mathbf{z}}$} \centerline{and $\gamma^1_{f_t,\mathbf{z}}$ constant}
\end{mybox}}
\ar@{=>}@<22pt>[r]
& 
\mbox{
\begin{mybox}[colback=white,width = 4.5cm]{\textbf{\mbox{\eqref{cmf}}}}
\centerline{Diffeomorphism type of} \centerline{Milnor fibration of $f_t$ constant}
\end{mybox}}
\ar@{<=}@<22pt>[r]
\ar@{<=}[u]
& 
\mbox{
\begin{mybox}[colback=white,width = 3.2cm]{\textbf{\mbox{\eqref{c3p}}}}
\centerline{$\mathring\mu({f_t})$ and $\mu({f_t}+z_1^j)$}
\centerline{constant} 
\end{mybox}}
\ar@{<=>}[u]
}
\end{eqnarray*}
}
\caption{Numerical control of the topological type in a family of line singularities}
\label{diagram}
\end{figure}

Here, the word ``constant'' means ``independent of $t$ for all small $t$''. Conditions \eqref{c0}--\eqref{equimult} are related to each other as described in items \ref{fp}--\ref{lp} below (see also Figure~\ref{diagram} for a full overview).

\begin{property}\label{fp}
If $n\geq 5$, then 
$\eqref{c0}\Rightarrow \eqref{c1}$ and $\eqref{c0}\Rightarrow \eqref{c1p}$.
\end{property}

The implication $\eqref{c0}\Rightarrow \eqref{c1}$ is given in \cite[Proposition, p.~47]{M7} as a consequence of \cite[Theorem, p.~437]{Ti}. The relation $\eqref{c0}\Rightarrow \eqref{c1p}$ follows. Indeed, by $\eqref{c0}$, the Milnor number $\mu({f_t}\vert_{V(z_1)})$ is constant, and by \cite[Theorem 2.1]{LR}, this implies that the family $\{{f_t}\vert_{V(z_1)}\}$ is topologically  equisingular.

\begin{rem}
Here, and hereafter, the condition $n\geq 5$ comes from the use of the L\^e-Ramanujam and Timourian theorems (cf.~\cite{LR,Ti}), which themselves use the $h$-cobordism theorem (cf.~\cite{Smale}).
\end{rem}

\begin{property}\label{pa2}
$\eqref{c1}\Rightarrow \eqref{c8} \Leftrightarrow \eqref{c7}$; if $n\geq 5$, then $\eqref{c7}\Rightarrow \eqref{c1}$.
\end{property}

The implication $\eqref{c1}\Rightarrow \eqref{c8}$ is proved in \cite[Proposition, p.~380]{M5} and \cite[Sec.~4]{M7}. The equivalence $\eqref{c8}\Leftrightarrow \eqref{c7}$ is shown in \cite[Sec.~4]{M7}. The implication $\eqref{c7}\Rightarrow \eqref{c1}$ for $n\geq 5$ is proved in \cite[Theorem 42]{B}.

\begin{property}
If $n\geq 5$, then $\eqref{c0}\Rightarrow \eqref{c3}$
\end{property}

If $\eqref{c0}$ holds, then $\mathring\mu({f_t})$ and $\mu({f_t}\vert_{V(z_1)})$ are constant. Now, since $n\geq 5$ and since in this case $\eqref{c0}$ implies~$\eqref{c1}$ (which in turn implies $\eqref{c7}$), it follows from the uniform Iomdine-L\^e-Massey formula (see Pro\-position 2.1 and the relation (2.2) in \cite{M7} and Theorem 4.15 in \cite{M}) that for all $t$ sufficiently small and all $j$ sufficiently large, the function $f_t+z_1^j$ has an isolated singularity at $\mathbf{0}$ and its Milnor number, which is given by
\begin{equation*}
\mu(f_t+z_1^j) = \lambda^0_{f_t,\mathbf{z}}+(j-1)\lambda^1_{f_t,\mathbf{z}},
\end{equation*}
is constant.

\begin{property}\label{p14}
$\eqref{c7}\Leftrightarrow \eqref{c7bis}$
\end{property}

By \cite[Sec.~1]{M7}, $\lambda^1_{f_t,\mathbf{z}}=\mathring\mu({f_t}):=\mu({f_t}\vert_{V(z_1-a_1)})$ for any sufficiently small $a_1\not=0$.

\begin{property}\label{pa5}
$\eqref{c2}\Rightarrow \eqref{c4}$; if $n\geq 5$, then 
$\eqref{c2}\Leftrightarrow \eqref{c4}$.
\end{property}

Since $\eqref{c2}$ contains $\eqref{c1}$, it implies $\eqref{c7}$ (cf.~\ref{pa2}). As the Milnor number is a topological invariant, $\eqref{c2}$ also implies that $\mu({f_t}\vert_{V(z_1)})$ is constant. So, altogether, $\eqref{c2}\Rightarrow \eqref{c4}$. 

To show the converse when $n\geq 5$, first observe that in this case $\eqref{c4}$ (which contains $\eqref{c7}$) implies $\eqref{c1}$ (cf.~\ref{pa2}), which is the first half part of $\eqref{c2}$. To get the second part, remind that $\lambda^1_{f_t,\mathbf{z}}=\mathring\mu({f_t}):=\mu({f_t}\vert_{V(z_1-a_1)})$, where $a_1\not=0$ is small enough. Therefore, if $\lambda^1_{f_t,\mathbf{z}}$ is constant and $n\geq 5$, then, by \cite[Theorem~2.1]{LR}, for any sufficiently small $a_1\not=0$, the family $\{{f_t}\vert_{V(z_1-a_1)}\}$ is topologically  equisingular. That the family $\{{f_t}\vert_{V(z_1)}\}$ is topologically  equisingular too follows from the constancy of $\mu({f_t}\vert_{V(z_1)})$, the assumption $n\geq 5$ and \cite[Theorem 2.1]{LR} again. 

\begin{property}\label{a6}
$\eqref{c1p}\Rightarrow\eqref{c4}$; if $n\geq 5$, then $\eqref{c1p}\Leftrightarrow \eqref{c2}$
\end{property}

Since $\eqref{c1}\Rightarrow \eqref{c7}$ (cf.~\ref{pa2}) and the topological  equisingularity of $\{{f_t}\vert_{V(z_1)}\}$ implies the constancy of the Milnor number $\mu({f_t}\vert_{V(z_1)})$, the condition $\eqref{c1p}$ implies $\eqref{c4}$, which is equivalent to $\eqref{c2}$ if $n\geq 5$ (cf.~\ref{pa5}). 

\begin{property}\label{pc5ic6}
$\eqref{c5}\Leftrightarrow \eqref{c6}$
\end{property}

This is proved in \cite{M7}. Since the argument will be useful for us later, let us briefly recall it. By \cite[Proposition 1.23]{M}, $\gamma^1_{f_t,\mathbf{z}}+ \lambda^0_{f_t,\mathbf{z}}
= \bigl([\Gamma^1_{f_t,\mathbf{z}}] \cdot [V(f_t)]\bigr)_{\mathbf{0}}$.
By \cite[Corollary~2.4]{M7}, if \eqref{c5} holds, then for any integer $j$ sufficiently large, the Milnor numbers $\mu(f_t+z_1^j)$ and $\mu(f_t\vert_{V(z_1)})$ are constant. Indeed, by the uniform Iomdine-L\^e-Massey formula, for all $t$ sufficiently small and all $j$ sufficiently large, the function $f_t+z_1^j$ has an isolated singularity at $\mathbf{0}$ and we have:
\begin{align*}
& \tag{\ref{pc5ic6}.1}\label{fu1}
\mu(f_t+z_1^j) = \lambda^0_{f_t,\mathbf{z}}+(j-1)\lambda^1_{f_t,\mathbf{z}};
\\
&\mu(f_t+z_1^j)+\mu(f_t\vert_{V(z_1)}) =(\gamma^1_{f_t,\mathbf{z}}+\lambda^0_{f_t,\mathbf{z}})+ j \lambda^1_{f_t,\mathbf{z}}.\tag{\ref{pc5ic6}.2}\label{fu2}
\end{align*}
Thus, if \eqref{c5} holds, then the sum $\mu(f_t+z_1^j)+\mu(f_t\vert_{V(z_1)})$ is constant, and by the upper-semicontinuity of the Milnor number, this implies that both $\mu(f_t+z_1^j)$ and $\mu(f_t\vert_{V(z_1)})$ are constant. The condition \eqref{c6} follows immediately.

\begin{property}
$\eqref{c3}\Leftrightarrow \eqref{c4}$
\end{property}

If $\eqref{c3}$ holds, then $\lambda^1_{f_t,\mathbf{z}}=\mathring\mu({f_t})$ is constant. That $\lambda^0_{f_t,\mathbf{z}}$ is constant too follows from \eqref{fu1}. The converse $\eqref{c4}\Rightarrow \eqref{c3}$  follows  exactly from the same formula.

\begin{property}
$\eqref{c4}\Leftrightarrow \eqref{c5}$
\end{property}

If $\eqref{c4}$ holds, then \eqref{fu1} says that $\mu(f_t+z_1^j)$ is constant. Then, by \eqref{fu2}, $\gamma^1_{f_t,\mathbf{z}}+\lambda^0_{f_t,\mathbf{z}}$ is constant. Conversely, if $\eqref{c5}$ holds, then, by \eqref{fu2} and the upper-semicontinuity of the Milnor number, both $\mu(f_t+z_1^j)$ and $\mu(f_t\vert_{V(z_1)})$ are constant. Combined with \eqref{fu1}, this implies that $\lambda^0_{f_t,\mathbf{z}}$ is constant.

\begin{property}\label{pv}
$\eqref{c3p}\Leftrightarrow \eqref{c8}$ 
\end{property}

This is proved in \cite[Remarque 1, p.~544]{Va}.

\begin{rem}
In \cite[Remarque 3, p.~544]{Va}, Vannier gives an example of a $1$-dimen\-sional singularity (the singular set is a union of three lines intersecting at $\mathbf{0}$) for which  $\eqref{c3p}$ holds but $\eqref{c3}$ fails. We do not know any example satisfying $\eqref{c3p}$ and for which $\eqref{c3}$ does not hold when the singular set is just a line.
\end{rem}

\begin{property}
$\eqref{LG}\Leftrightarrow \eqref{c3}$ 
\end{property}

By the L\^e-Greuel formula \cite{L2,G2}, $m(f_t+z_1^j)=\mu({f_t}+z_1^j)+\mu({f_t}\vert_{V(z_1)})$. The result then follows from the upper-semicontinuity of the Milnor number.

\begin{property}
$\eqref{c3p}\Rightarrow \eqref{cmf}$ and $\eqref{c7}\Rightarrow \eqref{cmf}$
\end{property}

The implication $\eqref{c3p}\Rightarrow \eqref{cmf}$ is proved in \cite[Th\'eor\`eme B]{Va} while $\eqref{c7}\Rightarrow \eqref{cmf}$ is proved in \cite[Theorem 9.4]{M}. (Note  that $\eqref{c3p}\Leftrightarrow \eqref{c7}$ by \ref{pv} and \ref{pa2}.)

\begin{property}\label{lp}
$\eqref{c00}\Rightarrow \eqref{equimult}$; 
$\eqref{c00}\Rightarrow \eqref{c1}$;
$\eqref{c5}\not\Rightarrow \eqref{c00}$;
if $n=3$, then $\eqref{c00}\Leftrightarrow \eqref{cr1}$.
\end{property}

$\eqref{c00}\Rightarrow \eqref{equimult}$ is proved in \cite[Corollary 6.2]{H}; $\eqref{c00}\Rightarrow \eqref{c1}$ is an immediate consequence of the Thom-Mather first isotopy theorem \cite{Th,Ma}; $\eqref{c5}\not\Rightarrow \eqref{c00}$ is proved in \cite[Sec.~5]{M7}; finally, $\eqref{c00}\Leftrightarrow \eqref{cr1}$ for $n=3$ is proved in \cite[Corollary 6.6]{GG}.

\begin{rem}
Combined with \ref{pa2}, it follows that $\eqref{cr1}\Rightarrow \eqref{c1}$ if $n\not=4$.
\end{rem}

\bibliographystyle{amsplain}

\quad

\bigskip
\footnotesize
\noindent
\textbf{Christophe Eyral}\\
Institute of Mathematics,
Polish Academy of Sciences\\
\'Sniadeckich 8,
00-656 Warszawa, Poland\\
\textit{E-mail:} ch.eyral@impan.pl

\vskip 2mm

\noindent
\textbf{Maria Aparecida Soares Ruas}\\
Instituto de Ci\^encias Matem\'aticas e de Computa\c c\~ao,
Universidade de S\~ao Paulo\\
Avenida Trabalhador S\~ao-Carlense, 400 - Centro,
Caixa Postal 668\\ 
13566-590 S\~ao Carlos, SP, Brazil\\
\textit{E-mail:} maasruas@icmc.usp.br

\end{document}